\documentclass[11pt]{amsart}
\usepackage{amsfonts}
\usepackage{amsthm}

\usepackage{a4wide}

\usepackage{amsmath,amssymb,enumerate}
\usepackage{epsfig,fancyhdr,color}

\usepackage{amssymb}
\usepackage{amsmath,amsthm}
\usepackage{latexsym}
\usepackage{amscd}
\usepackage{psfrag}
\usepackage{graphicx}
\usepackage[latin1]{inputenc}
\usepackage[all]{xy}
\usepackage[mathcal]{eucal}

\definecolor{NoteColor}{rgb}{1,0,0}

\renewcommand{\textsc}{\textcolor{red}}

\newtheorem{theorem}{\rm\bf Theorem}[section]

\newtheorem*{theorem 1}{\rm\bf Proposition 1}
\newtheorem*{theorem 2}{\rm\bf Proposition 2}

\theoremstyle{definition}
\newtheorem{definition}[theorem]{\rm\bf Definition}

\theoremstyle{remark}

\newcommand{\mK}{{{\mathbb K}}}
\newcommand{\mH}{{{\mathbb H}}}
\newcommand{\mO}{{{\mathbb O}}}

\newcommand{\calL}{{{\mathcal L}}}
\newcommand{\calQ}{{{\mathcal Q}}}

\newcommand{\calS}{{{\mathcal S}}}

\newcommand{\C}{{{\mathbb C}}}
\newcommand{\R}{{{\mathbb R}}}
\newcommand{\mN}{{{\mathbb N}}}

\newcommand{\bH}{{{\bf H}}}
\newcommand{\bc}{{{\bf c}}}
\newcommand{\bp}{{{\bf p}}}

\newcommand{\calC}{{{\mathcal C}}}

\newcommand{\fn}{{{\mathfrak n}}}

\newcommand{\fh}{{{\mathfrak H}}}

\def\Barint_#1{\mathchoice
          {\mathop{\vrule width 6pt height 3 pt depth -2.5pt
                  \kern -8pt \intop}\nolimits_{#1}}%
          {\mathop{\vrule width 5pt height 3 pt depth -2.6pt
                  \kern -6pt \intop}\nolimits_{#1}}%
          {\mathop{\vrule width 5pt height 3 pt depth -2.6pt
                  \kern -6pt \intop}\nolimits_{#1}}%
          {\mathop{\vrule width 5pt height 3 pt depth -2.6pt
                  \kern -6pt \intop}\nolimits_{#1}}}

\begin{document}

\title{Quasiconformal mappings on the Heisenberg group: An overview}

\author{Ioannis D. Platis}

\address{
Department of Mathematics and Applied Mathematics, University of Crete\\
University Campus, Voutes,  70013 Heraklion Crete, Greece\\
email:\,\tt{jplatis@math.uoc.gr}
}

\begin{abstract} We present a brief overview of the Kor\'anyi-Reimann \index{Kor\'anyi-Reimann theory} theory of quasiconformal mappings  \index{quasiconformal mapping}\index{mapping!quasiconformal} on the Heisenberg group stressing on the analogies as well as on the differences between the Heisenberg group case and the classical two-dimensional case. We examine the extensions of the theory to more general spaces and we state some known results and open problems. 
\end{abstract}

\email{jplatis@math.uoc.gr}

\maketitle

\tableofcontents


 \section{Introduction}\label{s-1}
 The classical theory of quasiconformal mappings  \index{quasiconformal mapping}\index{mapping!quasiconformal} on the plane was developed first in the Euclidean space $\R^n$ and produced a variety of results, most of them closely connected to topics in Analysis. In the non-Riemannian setting, the foundations of the theory are tracked down in the works of Mostow \index{Mostow, G.D.} and Pansu. \index{Pansu, P.} In his celebrated 1968 rigidity result \cite{mos2}, Mostow \index{Mostow, G.D.} proved that in dimensions $n>2$ diffeomorphic compact Riemannian manifolds with constant negative curvature are isometric, in particular they are conformally equivalent. The proof of this result relies on the use of quasiconformal mappings  \index{quasiconformal mapping}\index{mapping!quasiconformal} of $\R^n$. Later, see \cite{M}, he extended this result to the setting of symmetric spaces of rank 1 of non-compact type, i.e., hyperbolic spaces $\bH^n_\mK$ where $\mK$ can be  the set of real numbers $\R$ (except when $n=2$), 
the set of complex numbers $\C$, the set of quaternions $\mH$, or the set of octonions $\mO$ (the latter only when $n=2$). To obtain this, Mostow \index{Mostow, G.D.} had to develop quasiconformal mappings  \index{quasiconformal mapping}\index{mapping!quasiconformal} on the boundary of these spaces as an indispensable tool. In rough lines, Mostow's \index{Mostow, G.D.} proof in the case $\mK=\C$ and $n=2$ 
goes 
as follows. Let $G$ and $G'$ be two cocompact lattices, i.e., $M=\bH^2_\C/G$ and $M'=\bH^2_\C/G'$ are compact and suppose that $\rho:G\to G'$ is an isomorphism. From $\rho$ it is possible to define 
a quasi-isometric self map $F$ of $\bH^2_\C$ which is equivariant; this map need not be even continuous but has the property that it takes geodesics to quasi-geodesics. Due to a fundamental result in Gromov \index{Gromov, M.} hyperbolic spaces,\index{hyperbolic space!Gromov} from this property $F$ is extended to a boundary map $F_\infty:\partial\bH^2_\C\to\bH^2_\C$ which is in fact a quasiconformal homeomorphism of $\partial\bH^2_\C$ with respect to a metric comparable with the Kor\'anyi-Cygan metric. \index{Kor\'anyi-Cygan metric}\index{metric!Kor\'anyi-Cygan} The latter is defined on the (first) Heisenberg group \index{Heisenberg group} $\fh^1_\C$; this group is the $N$ group in the $KAN$ decomposition of the symmetric space $\bH^2_\C$ and $\partial\bH^2_\C$ is the one point compactification of $\fh$. After showing that $F_\infty$ has enough regularity, Mostow \index{Mostow, G.D.} proves that a $(G,G')$-equivariant quasiconformal self map of the boundary is associated with the action of an element of the 
isometry group ${\rm SU}(2,1)$ of $\bH^2_\C={\rm SU}(2,1)/{\rm SU}(2)$ and the proof is concluded from the equivariance of the resulting isometry.
We refer the interested reader to \cite{cap}, pp. 135--140, for a short but more detailed description of Mostow's \index{Mostow, G.D.} proof. 
Pansu \index{Pansu, P.} obtained a stronger rigidity statement for the cases $\mK=\mH$ and $\mK=\mO$ in \cite{Pa}. By using Mostow's \index{Mostow, G.D.} methods, Pansu \index{Pansu, P.} proved the following property which does not
hold for real and complex hyperbolic spaces: every quasi-isometry of quaternionic or
octonionic spaces has bounded distance from an isometry. Using the conformal geometry of
the boundary which is modelled on the nilpotent group \index{nilpotent group}\index{group!nilpotent} $\fh^{n-1}_\mK$ in the $KAN$ decomposition of the symmetric space $\bH^n_\mK$ endowed with a Carnot-Carath\'eodory metric \index{Carnot-Carath\'eodory metric} \index{metric!Carnot-Carath\'eodory} (i.e., a Carnot group), \index{Carnot group}\index{group!Carnot} and general properties of Loewner spaces, \index{Loewner space} (for the definitions of Carnot groups \index{Carnot group}\index{group!Carnot} and Loewner spaces, see Section \ref{sec:further}), he proved that {\it any} quasiconformal (in fact quasisymmetric) \index{quasisymmetric mapping}\index{mapping!quasisymmetric} homeomorphism of $\fh^{n-1}_\mK$ is actually conformal. This result does not apply to the case $\mK=\C$ and it is an open problem to understand the intrinsic reason for this phenomenon.

Mostow's rigidity \index{Mostow's rigidity theorem} had serious consequences; perhaps two of them are the most significant: First, the moduli space 
of hyperbolic metrics on a surface $\Sigma$, i.e., the  Teichm\"uller space \index{Teichm\"uller space} $\mathcal{T}(\Sigma)$, which is the case $\mK=\R$ and $n=2$ in the above setting, is just a counterexample of rigidity. The proof fails there, since it involves absolute continuity in measure of the boundary quasiconformal (actually quasisymmetric) mappings: in $S^1=\partial\bH^2_\R$ this does not hold.  

Second, the theory of quasiconformal mappings  \index{quasiconformal mapping}\index{mapping!quasiconformal} on the Heisenberg group \index{Heisenberg group} emerged, after the pioneering articles of Kor\'anyi-Reimann, \index{Kor\'anyi-Reimann theory} \cite{KR1} and  \cite{KR2}, and  Pansu, \cite{Pa}. \index{Pansu, P.} These works constituted a complete framework for the theory of quasiconformal mappings  \index{quasiconformal mapping}\index{mapping!quasiconformal} on the Heisenberg group \index{Heisenberg group} $\fh^{n-1}:=\fh^{n-1}_\C$. Especially, the exposition of Kor\'anyi \index{Kor\'anyi, A.} and Reimann  \index{Reimann, H. M.} is strongly influenced by the then state-of-the-art concerning the quasiconformal mappings  \index{quasiconformal mapping}\index{mapping!quasiconformal} in $\R^n$. There are many occasions where quasiconformal mappings  \index{quasiconformal mapping}\index{mapping!quasiconformal} on the Heisenberg group \index{Heisenberg group} behave in the same way as those of $\R^n$. On the 
other hand, there are significant differences: For instance, in Euclidean spaces $(n>2)$ there is no Beltrami equation whereas a system of Beltrami equations appears in the Heisenberg group \index{Heisenberg group} case for all $n$. However, in contrast to the real two-dimensional case where there exist solutions to the Beltrami equation, the solvability of the Beltrami system \index{Beltrami equation} in the Heisenberg group \index{Heisenberg group} case is still an open problem. According to Kor\'anyi \index{Kor\'anyi, A.} and Reimann, \index{Reimann, H. M.} even in the case of a full solution of the Beltrami system \index{Beltrami equation} this is not likely to produce startling results similar to the ones of the plane case. Regularity issues for the quasiconformal maps on the Heisenberg group \index{Heisenberg group} were first overlooked by Mostow \index{Mostow, G.D.} in his original rigidity result. It turned out that the basic property of absolute continuity on lines was much more difficult to obtain 
than in the Euclidean case. Kor\'anyi \index{Kor\'anyi, A.} and Reimann \index{Reimann, H. M.} brought this into Mostow's attention and together they remedied the problem; this correction appears in \cite{KR2}. Crucial to the further development of the theory was Pansu's \index{Pansu, P.} differentiability theorem; as in the Euclidean case, this theorem is derived from the Rademacher-Stepanov Theorem. Pansu's \index{Pansu, P.} notion of differentiability turned out to be completely adapted to the structure on the Heisenberg group; \index{Heisenberg group} 
the Pansu  \index{Pansu, P.} derivative is founded in such a way that it preserves the grading of the Lie algebra. 

All the above ignited the research of quasiconformal mappings \index{quasiconformal mapping}\index{mapping!quasiconformal} in various other spaces; only some of them are Carnot groups, \index{Carnot group}\index{group!Carnot} sub-Riemannian manifolds, \index{sub-Riemannian manifold}\index{manifold!sub-Riemannian} metric spaces with controlled geometry, etc.

In what follows we will briefly survey the theory of quasiconformal mappings \index{quasiconformal mapping}\index{mapping!quasiconformal} of the first Heisenberg group \index{Heisenberg group} $\fh$, that is, $\fh^1_\C$.  In Section \ref{sec:heis} we define the Heisenberg group \index{Heisenberg group} $\fh$ and comment on  its structures and properties. As the title of Section \ref{sec:overview} indicates, we overview there some basic results of the Kor\'anyi-Reimann \index{Kor\'anyi-Reimann theory} theory. Finally, in Section \ref{sec:further} further developments of the theory are briefly described. We finally underline here that this last section by no means contains all the developments of the theory; it rather reflects to the interests of the author.
\section{Heisenberg group}\label{sec:heis}

The (first) Heisenberg group \index{Heisenberg group} $\fh$ is the set $\C\times\R$ with multiplication $*$ given by
$$
(z,t)*(w,s)=(z+w,t+s+2\Im(z\overline{w})),
$$
for every $(z,t)$ and $(w,s)$ in $\fh$. We consider two metrics defined on $\fh$; the first one is induced via the Kor\'anyi map $\alpha:\fh\to \C$ which is given for every $(z,t)\in\fh$ by
$$
\alpha(z,t)=-|z|^2+it.
$$
Now, the Kor\'anyi gauge \index{Kor\'anyi gauge} $|\cdot|_\fh$ is defined by 
$$
\left|(z,t)\right|_\fh=|\alpha(z,t)|^{1/2}=\left| |z|^2-it\right|^{1/2},
$$
for every $(z,t)\in\fh$. Then the  {\it Kor\'anyi-Cygan} (or {\it Heisenberg})  metric \index{Kor\'anyi-Cygan metric}\index{metric!Kor\'anyi-Cygan} $d_\fh$ is defined by the relation
$$
d_\fh\left((z_1,t_1),\,(z_2,t_2)\right)
=\left|(z_1,t_1)^{-1}*(z_2,t_2)\right|.
$$
Note that the $d_\fh$-sphere of radius $R>0$ and centered at the origin, called  the {\it Kor\'anyi sphere}, is
$$
\calS_\fh(R)=\{(z,t)\in\fh\;|\;\left|(z,t)\right|_\fh=R\}.
$$
The metric $d_\fh$ is invariant under:
\begin{enumerate}
 \item {\it Left translations}  $T_{(\zeta,s)}$, $(\zeta,s)\in\fh$, that is, left actions of $\fh$ onto itself, which are given by
 $$
 T_{(\zeta,s)}(z,t)=(\zeta,s)*(z,t),
 $$
 for every $(z,t)\in\fh$;
 \item {\it rotations} $R_\theta$, $\theta\in\R$, around the vertical axis $\mathcal{V}=\{0\}\times\R$, 
 which are given by 
 $$
 R_\theta(z,t)=(ze^{i\theta},t),
 $$
 for every $(z,t)\in\fh$.
\end{enumerate}
Left translations are left actions of $\fh$ onto itself and rotations are induced by an action of ${\rm U}(1)$ on $\fh$; together they form the group ${\rm Isom}^+(\fh,d_\fh)$ of {\it orientation-preserving Heisenberg isometries}. The full group  ${\rm Isom}(\fh,d_\fh)$ of {\it  Heisenberg isometries} comprises compositions of elements  of  ${\rm Isom}^+(\fh,d_\fh)$ with the 
\begin{enumerate}
 \item [{3)}] {\it conjugation} $J$ which is defined by $$J(z,t)=(\overline{z},-t),$$
 for every $(z,t)\in\fh$.
\end{enumerate}

We also consider two other kinds of transformations, namely: 
\begin{enumerate}
 \item [{4)}] {\it Dilations} $ D_\delta$, $\delta>0$. These are defined by
 $$
 D_\delta(z,t)=(\delta z, \delta^2 t),
 $$
for every $(z,t)\in\fh$;
\item [{5)}] {\it inversion} $I$, which is defined in $\fh\setminus\{(0,0)\}$ by 
$$
I(z,t)=\left(z\left(\alpha(z,t)\right)^{-1},
-t\left|\alpha(z,t)\right|_\fh^{-2}\right).
$$
\end{enumerate}
The metric $d_\fh$ is  scaled up to multiplicative constants by the action of  dilations. Also, for each $p,q\in\fh\setminus\{o=(0,0)\}$ we have
$$
d_\fh(I(p),I(q))=\frac{d_\fh(p,q)}{d_\fh(o,p)\cdot d_\fh(o,q)}.
$$
Compositions of orientation-preserving Heisenberg isometries, dilations and inversion form the {\it similarity group} ${\rm Sim}(\fh)$ of $\fh$. Recall the {\it complex hyperbolic plane} \index{complex hyperbolic plane} $\bH^2_\C$: That is, the symmetric space ${\rm SU}(2,1)/{\rm SU}(2)$. The Heisenberg group \index{Heisenberg group} $\fh$ is the $N$ group in its $KAN$ decomposition and the boundary $\partial\bH^2_\C$ is the one point compactification of $\fh$, that is, $\partial\bH^2_\C=\fh\cup\{\infty\}$. It can be proved that the action of ${\rm SU}(2,1)={\rm Isom}(\bH^2_\C)$ on the boundary $\partial\bH^2_\C=\fh\cup\{\infty\}$ is completely described by the transformations (1)-(5), that is, if $g$ is an isometry, then it is the composition of transformations of the form (1)-(5).

The Heisenberg group \index{Heisenberg group} $\fh$ is a two-step nilpotent Lie group with underlying manifold $\C\times\R$; its left translations are  
$T_{(\zeta,s)}$, $(\zeta,s)\in\fh$, as in (1). Let $(z=x+iy,t)$ be the coordinates of $\fh$ and consider the left invariant vector fields
\begin{eqnarray*}
X=\frac{\partial}{\partial x}+2y\frac{\partial}{\partial t},\quad Y=\frac{\partial}{\partial y}-2x\frac{\partial}{\partial t},
\quad T=\frac{\partial}{\partial t}
\end{eqnarray*}
and the complex fields
\begin{eqnarray*}
Z=\frac{1}{2}(X-i Y)=\frac{\partial}{\partial z}+i\overline{z}\frac{\partial}{\partial t}\quad
\overline{Z}=\frac{1}{2}(X+i Y)=\frac{\partial}{\partial\overline{z}}-iz\frac{\partial}{\partial t}.
\end{eqnarray*}
The Lie algebra of left invariant vector fields of $\fh$ has a grading $\mathfrak{h} = \mathfrak{v}_1\oplus \mathfrak{v}_2$ with
\begin{displaymath}
\mathfrak{v}_1 = \mathrm{span}_{\R}\{X, Y\}\quad \text{and}\quad \mathfrak{v}_2=\mathrm{span}_{\R}\{T\}.
\end{displaymath}
The only non-zero bracket relation is
$$
[X,Y]=4T.
$$
The Heisenberg group \index{Heisenberg group} constitutes the prototype of both contact and sub-Riemannian geometry. \index{sub-Riemannian geometry}\index{geometry!sub-Riemannian}
The contact form \index{contact form}\index{form!contact} $\omega$ of $\fh$ is defined as the unique 1-form satisfying $X,Y\in{\rm ker}\omega$, $\omega(T)=1$. Uniqueness here is modulo change of coordinates as it follows by the contact version of Darboux's Theorem. \index{Darboux's Theorem} The distribution in $\fh$ defined by the first layer $\mathfrak{v}_1:={\rm H}(\fh)$ is called the {\it horizontal distribution}. \index{horizontal distribution}\index{distribution!horizontal} 
Explicitly, in the Heisenberg coordinates $z=x+iy,t,$ we have
\begin{eqnarray*}
\omega=dt+2(xdy-ydx)=dt+2i\Im(\overline{z}dz).
\end{eqnarray*}
The {\it sub-Riemannian metric  \index{sub-Riemannian metric}\index{metric!sub-Riemannian}} $\langle\cdot,\cdot\rangle$ is defined on ${\rm H}(\fh)$ by the relations
$$
\langle X,X\rangle=\langle Y,Y\rangle=1,\quad \langle X,Y\rangle=\langle Y,X\rangle=0. 
$$
The corresponding norm shall be denoted by $\|\cdot\|$. The {\it Legendrian foliation} of $\fh$ is the foliation of $\fh$ by horizontal curves. \index{horizontal curve}\index{curve!horizontal} An absolutely continuous curve $\gamma:[a,b]\to \fh$ (in the Euclidean sense) with 
\begin{displaymath}
\gamma(\tau)=(\gamma_h(\tau),\gamma_3(\tau))\in\mathbb{C}\times \mathbb{R},
\end{displaymath}
 is called {\it horizontal} if $\dot{\gamma}(\tau)\in {\rm H}_{\gamma(\tau)}(\fh)$ for almost all $\tau\in [a,b]$; equivalently,
 \begin{displaymath}
 \dot t(\tau)=-2\Im\left(\overline{z(\tau)}\dot z(\tau)\right),
\end{displaymath}
for almost all $\tau\in [a,b]$.
A curve $\gamma:[a,b]\to \fh$ is absolutely continuous with respect to $d_\fh$ if and only if it is a horizontal curve. \index{horizontal curve}\index{curve!horizontal} Moreover,
the horizontal length of a smooth rectifiable curve $\gamma=(\gamma_h,\gamma_3)$ with respect to $\|\cdot\|$ is given by the integral over the (Euclidean) norm of the horizontal part of the tangent vector,
\begin{displaymath}
 \ell_h(\gamma)=\int_a^b \|\dot{\gamma}_h(\tau)\|\;d\tau=\int_a^b\left(\langle\dot\gamma(\tau),X_{\gamma(\tau)}\rangle^2+\langle\dot\gamma(\tau),Y_{\gamma(\tau)}\rangle^2\right)^{1/2}d\tau.
\end{displaymath}
Thus the sub-Riemannian or  {\it Carnot-Carath\'eodory distance} \index{Carnot-Carath\'eodory distance}\index{distance!Carnot-Carath\'eodory} of two arbitrary points $p,q\in\fh$ is
$$
d_{cc}(p,q)=\inf_\gamma\ell_h(\gamma),
$$ 
where $\gamma$ is horizontal and joins $p$ and $q$ (and therefore horizontal curves \index{horizontal curve}\index{curve!horizontal} are {\it geodesics} for the sub-Riemannian metric). \index{sub-Riemannian metric}\index{metric!sub-Riemannian} The $d_{cc}$-sphere \index{Carnot-Carath\'eodory sphere}\index{sphere!Carnot-Carath\'eodory} of radius $R$ and centered at the origin is the {\it Carnot-Carath\'eodory sphere} $\calS_{cc}(R)$ and is constructed as follows. Consider the family ${\bf c}_k$, $k\in\R$, of planar circles where
$$
{\bf c}_k(s)=\frac{1}{k}(1-e^{iks}),\quad s\in[0,2\pi/|k|].
$$
In the case where $k=0$, circles degenerate into straight line segments. For every $k$, ${\bf c}_k$ is lifted to the horizontal curve \index{horizontal curve}\index{curve!horizontal} $\bp_0^k$ where
$$
\bp_0^k(s)=\left(\bc_k(s),t_k(s)\right),\quad t_k(s)=\frac{2}{k}\left(\frac{1}{k}\sin(ks)-s 
\right).
$$
Denote by $\bp_\phi^k$ the rotation of $\bp_0^k$ around the vertical axis $\mathcal{V}=\{(z,t)\in\fh\;|\;z=0\}$. Then
$$
\calS_{cc}(R)=\bp_\phi^k(R),\;k\in[-2\pi/R,2\pi/R],\;\phi\in[0,2\pi].
$$
In comparison, the Kor\'anyi-Cygan metric \index{Kor\'anyi-Cygan metric}\index{metric!Kor\'anyi-Cygan} and the Carnot-Carath\'eodory metric \index{Carnot-Carath\'eodory metric} \index{metric!Carnot-Carath\'eodory} are equivalent metrics, they both behave like the Euclidean metric in horizontal directions ($X$ and $Y$), and behave like the square root of the Euclidean metric in the missing direction $(T)$. Their isometry groups are the same, but not their similarity groups: inversion is not a similarity of $d_{cc}$.
The relation of $d_\fh$ (which is not a geodesic metric) and $d_{cc}$ is given as follows, see \cite{cap}: there exist constants $C_1,C_2>0$ so that
\begin{equation}\label{eq:metrics}
C_1d_\fh(p,0)\le d_{cc}(p,0)\le C_2d_\fh(p,0),
\end{equation}
for each $p\in\fh$. Finally, $d_\fh$ and $d_{cc}$ generate the same infinitesimal structure: If $\gamma:[0,1]\to\fh$ is a $\calC^1$ curve and $t_i=i/n$, $i=1,\dots,n$, is a partition of $[0,1]$, then
$$
\limsup_{n\to\infty}\sum_{i=1}^n d_\fh\left(\gamma(t_i),\gamma(t_{i-1}\right)=
\left\{\begin{matrix}
        \ell(\gamma)&\text{if}\;\gamma\;\text{is horizontal},\\
        \infty&\text{otherwise}.
       \end{matrix}\right.
$$
Here, $\ell(\gamma)$ denotes the length of $\gamma$ with respect to both $d_\fh$ and $d_{cc}$, see Lemma 2.4 in \cite{cap}.

Contact transformations \index{contact mapping}\index{mapping!contact} between domains (open and connected subsets) of $\fh$ play an important role in the theory of quasiconformal mappings \index{quasiconformal mapping}\index{mapping!quasiconformal} of $\fh$. A {\it contact transformation} $f:\Omega\to\Omega'$ on $\fh$ is a diffeomorphism between domains  $\Omega$ and $\Omega'$ in $\fh$  which preserves the contact structure, i.e.,
\begin{equation}\label{eq:contact_form}
 f^*\omega = \lambda \omega,
\end{equation}
for some non-vanishing real valued function $\lambda$. We  write $f=(f_I,f_3)$, $f_I=f_1+\mathrm{i}f_2$. Then a contact mapping \index{contact mapping}\index{mapping!contact} $f$ is completely determined by $f_I$ in the sense that the contact condition \index{contact condition} (\ref{eq:contact_form}) is equivalent to the following system of partial differential equations:
\begin{eqnarray}
\label{eq:C1} 
&&\overline{f}_IZ f_I - f_I Z \overline{f}_I +iZf_3=0,\\
&&
\label{eq:C2}
f_I \overline{Z}\overline{f}_I-\overline{f}_I \overline{Z}f_I -
i\overline{Z}f_3=0,\\
&&
\label{eq:C3}
-i(\overline{f}_I T f_I - f_I T \overline{f}_I +
iTf_3)=\lambda.
\end{eqnarray}
If $f$ is $\calC^2$ it is elementary to prove that
$$
\det J_f=\lambda^2,
$$
where by $J_f$ we denote the usual Jacobian matrix  of $f$.

\section{A Brief Overview of the Kor\'anyi-Reimann Theory}\label{sec:overview}

We start by recalling various definitions of quasiconformal mappings  \index{quasiconformal mapping}\index{mapping!quasiconformal} on the Heisenberg group \index{Heisenberg group}\index{Heisenberg group!quasiconformal mappings} $\fh$. All these definitions turn out to be equivalent. In our exposition, we follow the lines of \cite{KR1} and \cite{KR2} with minor deviations.
\subsubsection*{Metric definition.}
 Kor\'anyi \index{Kor\'anyi, A.} and Reimann \index{Reimann, H. M.} wrote in \cite{KR1}: {\it The (metric) definition is a straightforward generalisation of the corresponding notion in Euclidean space. It was Mostow \index{Mostow, G.D.} who for the first time studied these mappings in the context of semisimple Lie groups of rank one}. We consider the Heisenberg group $\fh$ equipped with its metric $d_\fh$. For a homeomorphism $f$ between domains $\Omega$ and $\Omega'$ in $\fh$ and for $p\in\Omega$ we set
\[
 L_f(p,r)=\sup_{d_\fh(p,q)\le r,\;q\in\Omega}d_\fh(f(p),f(q)),\]

\[ l_f(p,r)=\inf_{d_\fh(p,q)\ge r,\;q\in\Omega}d_\fh(f(p),f(q)),\]

\[ H_f(p,r)=\frac{L(p,r)}{l(p,r)}\]
and 
\[ H_f(p)=\limsup_{r\to 0}H(p,r).\]
\begin{definition}\label{definition:metric}{\bf (Metric definition)}
A homeomorphism $f:\Omega\to\Omega'$ between domains $\Omega$ and $\Omega'$ in $\fh$ is quasiconformal if it is uniformly bounded from above in $\Omega$. If in addition
$$
{\rm ess\; sup}_{p\in\Omega}H_f(p)=\|H_f\|_\infty\le K,
$$
then $f$ is called $K$-{\it quasiconformal}. 
\end{definition}

\medskip

We note that due to (\ref{eq:metrics}), we obtain an equivalent definition if we substitute $d_\fh$ by the metric $d_{cc}$, see for instance the metric definition \index{metric definition of quasiconformality} in \cite{KR2}. Moreover, conformal mappings (i.e., elements of ${\rm SU}(2,1)$, the isometry group of complex hyperbolic plane $\bH^2_\C$, acting on $\fh$) are 1-quasiconformal. The converse is also true; a 1-quasiconformal mapping is necessarily an element of ${\rm SU}(2,1)$, see \cite{KR1} for a proof.

In the context of the Heisenberg group, there are various equivalent analytic definitions \index{analytic definition of quasiconformality} of quasiconformal mappings  \index{quasiconformal mapping}\index{mapping!quasiconformal} $f:\Omega\to\fh$ which are all equivalent to the metric definition \index{metric definition of quasiconformality} \ref{definition:metric}. We refer the reader to \cite{D}, \cite{H1} and \cite{V}.
For the definition we are about to state, see \cite{BFP1} and the references therein.

\subsubsection*{Analytic definition.}
We first define $HW^{1,4}(\Omega, \fh)$, the {\it horizontal Sobolev space}. \index{Sobolev space!horizontal} We say that a function $u:\Omega\to\C$ is in $HW^{1,4}(\Omega, \fh)$ if it is in $L^{4}(\Omega, \fh)$ and if there exist functions $v,w\in L^{4}(\Omega, \fh)$ such that
$$
\int_\Omega v\phi d\calL^3=-\int_\Omega uZ\phi d\calL^3 \quad\text{and}\quad \int_\Omega w\phi d\calL^3=-\int_\Omega u\overline{Z}\phi d\calL^3 
$$
for all $\phi\in\mathcal{C}^\infty_0(\Omega, \C)$. Now a mapping  $f:\Omega\to\fh$, $f=(f_I,f_3)$ is said to be in $HW^{1,4}(\Omega, \fh)$ if both $f_I,f_3$ are in  $HW^{1,4}(\Omega, \fh)$.  Such a mapping which also satisfies Conditions (\ref{eq:C1}) and (\ref{eq:C2}) a.e. is called {\it weakly contact} \index{weakly contact mapping}\index{mapping!weakly contact} and one can define its formal horizontal differential $(D_h)f_p$ at almost all $p$, which in matrix form is given by
\begin{equation*}
 (D_h)f_p=\left(\begin{matrix} Zf_I & \overline{Z}f_I\\
   Z\overline{f_I}&\overline{Zf_I}\\
                \end{matrix}\right)_p.
\end{equation*}
This is extended to a Lie algebra homomorphism known as the {\it Pansu derivative} $(D_0)f_p$, see \cite{Pa},  which in matrix form is
\begin{equation*}
 (D_0)f_p=\left(\begin{matrix} Zf_I & \overline{Z}f_I& 0\\
                 Z\overline{f_I}&\overline{Zf_I}&0\\
                 0&0&|Zf_I|^2-|\overline{Z}f_I|^2
                \end{matrix}\right)_p.
\end{equation*}
It is worth noticing that Pansu \index{Pansu, P.} derivation \index{Pansu derivative} is the Heisenberg analogue of the plain  derivation in Euclidean spaces; a mapping $f:\Omega\to\Omega'$ between domains of $\fh$ is  $P$-differentiable at $p\in\Omega$ if for $c\to 0$ the mappings
$$
D_c^{-1}\circ T_{f(p)}^{-1}\circ f\circ T_{f(p)}\circ D_c
$$
converge locally uniformly to a homomorphism $(D_0)f_p$ from $T_p(\fh)$ to $T_{f(p)}(\fh)$  which preserves the horizontal space  \index{horizontal space}\index{space!horizontal} ${\rm H}(\fh)$. Here $D$ and $T$ are dilations and left translations respectively.

Now let
\begin{enumerate}
 \item $\|(D_h)f_p\|:=\max\left\{\|(D_h)f_p(V)\|\;|\; \|V\|=1\right\}=|Zf_I(p)|+|\overline{Z}f_I(p)|$ a.e.;
 \item $J_f(p)=\det (D_0)f_p=(\det (D_h)f_p)^2=\left(|Zf_I(p)|^2-|\overline{Z}f_I(p)|^2\right)^2$;
 \item $$
 K(p,f)^2=\frac{\|(D_h)f_p\|^4}{J_f(p)}=\left(\frac{|Zf_I(p)|+|\overline{Z}f_I(p)|}{|Zf_I(p)|-|\overline{Z}f_I(p)|}\right)^2.
 $$
\end{enumerate}
\begin{definition}\label{definition:analytic}
{\bf (Analytic definition)} A homeomorphism $f:\Omega\to\Omega'$, $f=(f_I,f_3)$, between domains in $\fh$ is an orientation-preserving $K$-quasiconformal mapping if $f\in HW^{1,4}(\Omega, \fh)$ is weakly contact and if 
\begin{equation*}
\|(D_h)f_p\|^4\le KJ_f(p),
\end{equation*}
for almost all $p$.
\end{definition}
The function $\Omega\ni p\mapsto K(p,f)\in [1,\infty)$ is the {\it distortion function \index{distortion function}\index{function!distortion}} of $f$ and the constant of quasiconformality $K$ is also called the {\it maximal distortion} of $f$.

We remark that the analytic definition \index{analytic definition of quasiconformality} given by Kor\'anyi \index{Kor\'anyi, A.} and Reimann,  \index{Reimann, H. M.} see for instance \cite{KR2}, is based on the notion of absolute continuity on lines (ACL):\index{absolute continuity in lines}  Mappings with this property are absolutely continuous on a.e. fiber of any smooth fibration determined by a horizontal vector field $V$.
In the case of the Heisenberg group $\fh$, absolute continuity holds on almost all fibers of smooth {\it horizontal} fibrations. For such a fibration, the fibers $\gamma_p$ can be parametrised by the flow $f_s$ of a horizontal {\it unit} vector field $V$: i.e., $V$ is of the form $aX+bY$ with $a^2+b^2=1$. Mostow \index{Mostow, G.D.} proved (Theorem A in \cite{KR2}) that  quasiconformal mappings  \index{quasiconformal mapping}\index{mapping!quasiconformal} are absolutely continuous on a.e. fiber $\gamma$ of any given fibration $\Gamma_V$ determined by a left invariant horizontal vector field $V$. The  Euclidean counterpart was proved by Gehring in \cite{Gehring1960}. The proof in the Euclidean case is considerably easier. In this way, Kor\'anyi \index{Kor\'anyi, A.} and Reimann \index{Reimann, H. M.} defined a homeomorphism $f:\Omega\to\Omega'$, $f=(f_I,f_3)$, between domains in $\fh$ to be an orientation-preserving $K$-quasiconformal mapping if
\begin{enumerate}
\item[{(i)}] it is ACL;
\item[{(ii)}] it is a.e. $P$-differentiable, and
\item[{(iii)}] it satisfies a.e. the system of Beltrami equations \index{Beltrami equation} 
\begin{eqnarray}\label{eq:B1}
&&
\overline{Z}f_I=\mu Zf_I,\\
&&\label{eq:B2}
\overline{Z}f_{II}=\mu Zf_{II},
\end{eqnarray}
where $f_{II}=f_3+i|f_I|^2$ and $\mu$ is a complex function in $\Omega$ such that 
$$
\frac{1+\|\mu\|_\infty}{1-\|\mu\|_\infty}\le K \quad \text{a.e.}
$$
Here, $\|\mu\|_\infty={\rm ess sup}\{|\mu(z,t)|\;|\;(z,t)\in\Omega\}$. For each $p=(z,t)\in\Omega$, the function
\begin{equation*}
 \mu(p)=\mu_f(p)=\frac{\overline{Z}f_I(p)}{Zf_I(p)}
\end{equation*}
is called the {\it Beltrami coefficient of $f$}.
\end{enumerate}

\subsubsection*{Geometric definition.}
Kor\'anyi \index{Kor\'anyi, A.} and Reimann  \index{Reimann, H. M.} proved the equivalence of the metric and analytic definitions by showing that both are equivalent to a third definition which they called the geometric definition. \index{geometric definition of quasiconformality} This involves the notion of {\it capacity} of a ring domain ${\rm Cap}(E,F)$ in $\fh$: by $(E,G)$ we denote the open bounded subset $U=\fh\setminus(E\cup F)$ where $E$ and $F$ are disjoint connected closed subsets of $\fh$ and moreover $E$ is compact. Then the capacity  \index{capacity} of $(E,G)$ is
$$
{\rm Cap}(E,F)=\inf\int_\fh|\nabla_h u|^4d\calL^3,
$$
where the infimum is taken over all smooth functions $u$ in $G$ with $u_{|E}\ge 1$ and $u_{|F}=0$. Here $\nabla_h u$ is the horizontal gradient $(Xu)X+(Yu)Y$, and $d\calL^3$ is the volume element of the usual Lebesgue measure in $\C\times\R$. For details, see Section 3 in \cite{KR2}. 

\begin{definition}\label{definition:geometric1}{\bf (Geometric definition I)}
Let $f:\Omega\to\Omega'$ be a homeomorphism between domains in $\fh$. Then $f$ is quasiconformal if there exists a $K\ge 1$ such that for each ring domain $(E,F)$ in $\Omega$ we have the capacity  inequality \index{capacity inequality}
\begin{equation*}
 {\rm Cap}(E,F)\le K^2{\rm Cap}f(E,F).
\end{equation*}
\end{definition}
In what follows we will give another geometric definition \index{geometric definition of quasiconformality} which involves the notion of {\it modulus} of families of curves. \index{modulus of family of curves} Let $\Gamma$ be a family %
of rectifiable (with respect to $d_\fh$) curves lying in a domain $\Omega\subset\fh$, 
i.e., the curves have finite length with respect to $d_\fh$ (or in other words, they have finite horizontal length). If $\rho:\fh\to[0,\infty)$ is a non-negative Borel function and $\gamma$ is a parametrisation of a rectifiable curve $\gamma(t)=(\gamma_h(t),\gamma_3(t))$, $t\in[a,b]$, we define
$$
\int_\gamma\rho ds=\int_a^b\rho(\gamma(t))\|\dot\gamma_h(t)\|dt.
$$
Let ${\rm adm}(\Gamma)$ be the set of these non-negative Borel functions $\rho$ defined in $\fh$ which satisfy
$$
\int_\gamma \rho ds\ge 1, \quad\text{for all rectifiable}\;\gamma\in\Gamma.
$$
Then the {\it modulus} \index{modulus of families of curves} ${\rm Mod}(\Gamma)$ of $\Gamma$ is defined as 
$$
 {\rm Mod}(\Gamma)=\inf_{\rho\in{\rm adm}(\Gamma)}\int_\fh \rho^4(p) d\calL^3(p),
$$
where $d\calL^3$ is the volume element of the usual Lebesgue measure in $\C\times\R$. 

A fundamental inequality which may be found for instance in \cite{BFP1}, is formulated as follows.
Let $f:\Omega\to\Omega'$ be a $K$-quasiconformal mapping between domains in $\fh$ and let $\Gamma$ be a family of curves in $\Omega$. Then in the first place we have
\begin{equation}\label{eq:mod-int}
 {\rm Mod}(f(\Gamma))\le\int_\Omega K(p,f)^2\rho^4(p) d\calL^3(p).
\end{equation}
Here, $K(p,f)$ is the distortion function \index{distortion function}\index{function!distortion} of $f$ and the integral on the right is called the {\it distortion functional \index{distortion functional}\index{functional!distortion}} of $f$.
From (\ref{eq:mod-int}) we then obtain
\begin{equation}\label{eq:modulus-ineq}
 \frac{1}{K^2}{\rm Mod}(\Gamma)\le {\rm Mod}(f(\Gamma))\le K^2{\rm Mod}(\Gamma).
 \end{equation}
Inequality (\ref{eq:modulus-ineq}) is known as the {\it modulus inequality}. \index{modulus inquality} As a direct corollary of the modulus inequality  \index{modulus inquality} we obtain that the modulus \index{modulus of families of curves} of a family of curves is a conformal invariant; if $f$ is conformal then ${\rm Mod}(f(\Gamma))={\rm Mod}(\Gamma)$. Now the second geometric definition \index{geometric definition of quasiconformality} of quasiconformality stands as follows.

\begin{definition}\label{definition:geometric2}{\bf (Geometric definition II)}
Let $f:\Omega\to\Omega'$ be a homeomorphism between domains in $\fh$. Then $f$ is quasiconformal if there exists a $K\ge 1$ such that for each curve family $\Gamma$ in $\Omega$ the modulus inequality  \index{modulus inquality} (\ref{eq:modulus-ineq}) holds.
\end{definition}

Since the capacity of a ring domain  \index{capacity of a ring domain} $(E,G)$ is equal to the modulus \index{modulus of families of curves} of the family of horizontal curves \index{horizontal curve}\index{curve!horizontal} joining $E$ and $F$ in $U$, see \cite{E} or Proposition 2.4 of \cite{HK1}, the geometric definition II implies the geometric definition I. The converse may be derived via {\it quasisymmetric mappings}. \index{quasisymmetric mapping}\index{mapping!quasisymmetric} Like in the classical case, quasiconformal mappings  \index{quasiconformal mapping}\index{mapping!quasiconformal} are strongly related to quasisymmetric ones. A mapping $f:\Omega\to\Omega'$ between domains of $\fh$ is called {\it locally $\eta$-quasisymmetric} if there exists an increasing self-homeomorphism $\eta$ of $[0,\infty)$ such that for each Whitney ball $B\subset\Omega$,
$$
\frac{d_\fh(f(p),f(q))}{d_\fh(f(p),f(r))}\le \eta\left(\frac{d_\fh(p,q)}{d_\fh(p,r)}\right),
$$
for all $p,q,r\in B$, $p\neq r$. Recall that a Whitney ball $B\subset\Omega$ is a closed metric ball $B(x,R)$ centered at $x\in\Omega$ and with radius $R$, such that $2B=B(x,2R)\subset\Omega$.

The equivalence of all the afore stated definitions of quasiconformality is clarified by a theorem which in its full generality is Theorem 9.8 in \cite{HKST}, see also Theorem 6.33 in \cite{cap}. This result states that if $f:\Omega\to\Omega'$ is a homeomorphism between domains of $\fh$, then the following are equivalent.
 \begin{enumerate}
  \item [{i)}] $f$ is quasiconformal according to the metric definition \ref{definition:metric};
  \item [{ii)}] $f$ is locally $\eta$-quasisymmetric;
  \item [{iii)}] $f$ is quasiconformal according to the geometric definition \ref{definition:geometric2}.
 \end{enumerate}

In this manner, it follows that i), ii) and iii) are all equivalent to the analytic definition \index{analytic definition of quasiconformality} \ref{definition:analytic} and to geometric definition \index{geometric definition of quasiconformality} \ref{definition:geometric1} as well. 

\medskip

\subsubsection*{Smooth quasiconformal mappings.} In their first paper \cite{KR1}, Kor\'anyi \index{Kor\'anyi, A.} and Reimann  \index{Reimann, H. M.} studied mostly smooth quasiconformal mappings. \index{quasiconformal mapping!smooth}  \index{quasiconformal mapping}\index{mapping!quasiconformal}
Quasiconformal mappings with sufficient smoothness have to be contact transformations. \index{contact mapping}\index{mapping!contact} This property distinguishes quasiconformal mappings  \index{quasiconformal mapping}\index{mapping!quasiconformal} on the Heisenberg group $\fh$ from those defined on Euclidean spaces. In fact, from $P$-differentiability of quasiconformal mappings  \index{quasiconformal mapping}\index{mapping!quasiconformal} it follows that $P$-diffeomorphic $K$-quasiconformal mappings  \index{quasiconformal mapping}\index{mapping!quasiconformal} are contact transformations satisfying
\begin{equation}\label{cond:dil}
 \|(D_h)f\|^4\le K |J_f|\quad \text{a.e.}
\end{equation}
(Here, the absolute value in the Jacobian covers both situations of orientation--preserving and orientation-reversing mappings). The converse is also true, see Proposition 8 in \cite{KR2}: if a $\calC^2$ contact transformation $f$ satisfies Condition (\ref{cond:dil}), then $f$ is $K$-quasiconformal. We conclude that $K$-quasiconformal diffeomorphisms lie in the class of contact transformations. Due to the contact Conditions \index{contact condition} (\ref{eq:C1}), (\ref{eq:C2}) and (\ref{eq:C3}), Equation (\ref{eq:B1}) in the Beltrami system \index{Beltrami equation} implies Equation (\ref{eq:B2}).

\subsubsection*{Quasiconformal deformations and measurable Riemann Mapping Theorem.} Perhaps one of the most striking results in the original work of Kor\'anyi \index{Kor\'anyi, A.} and Reimann, \index{Reimann, H. M.} is their generalisation to the Heisenberg setting of the famous measurable Riemann Mapping Theorem  \index{measurable Riemann Mapping Theorem} in its infinitesimal version. 
Besides its genuine importance, this theorem enables us to construct as many quasiconformal mappings \index{quasiconformal mapping}\index{mapping!quasiconformal} on $\fh$ as we wish, out of {\it quasiconformal deformations}. \index{quasiconformal deformation}Let $f_s:\fh\to\fh$, $f_s=f_s(z,t)$, $s\in\R$, be a $\calC^1$ one-parameter group of  transformations of $\fh$ with infinitesimal generator $V$, satisfying the initial condition $f_0(z,t)={\rm id}$. Then the following differential equation holds:
$$
\frac{d}{ds}f_s(z,t)=V(f_s(z,t)).
$$
We are interested primarily in one-parameter groups of {\it contact} transformations \index{contact mapping}\index{mapping!contact} since we  have seen that smooth enough quasiconformal mappings  \index{quasiconformal mapping}\index{mapping!quasiconformal} are contact. Infinitesimal generators of one-parameter groups of contact transformations have been studied by Libermann \index{Libermann, P.} and are of a special form: According to Theorem 5 in \cite{KR1},  
$C^1$ vector fields of the form
\begin{equation}\label{eq:V}
V=-\frac{1}{4}\left[(Yp)X-X(p)Y\right]+pT=\frac{i}{2}\left[(\overline{Z}p)Z-(Zp)\overline{Z}\right]+pT,
\end{equation}
where $p$ is an arbitrary real valued function, generate local one-parameter groups of contact transformations. Conversely, every $\calC^1$ vector field $V$ which generates a local one-parameter group of contact transformations is necessarily of this form with $p=\omega(V)$.

A precise estimate for the constant of quasiconformality of a one-parameter group of quasiconformal mappings  \index{quasiconformal mapping}\index{mapping!quasiconformal} generated by a $\calC^2$ vector field is given by Theorem 6 in \cite{KR1}; this is also the first version of the measurable Riemann mapping theorem for the Heisenberg group case:
Suppose that $V$ is a $\calC^2$ vector field of the form (\ref{eq:V}) generating a one-parameter group $\{f_s\}$ of contact transformations. If
\begin{equation*}
 |ZZp|\le c^2/2,
\end{equation*}
then $f_s$ is $K$-quasiconformal with the constant of quasiconformality $K=K(s)$ of $f_s$ satisfying
\begin{equation*}
 K+\frac{1}{K}\le 2 e^{c|s|}.
\end{equation*}
Kor\'anyi \index{Kor\'anyi, A.} and Reimann  \index{Reimann, H. M.} improved this result, see Theorem H in \cite{KR2}: The assumptions there for the vector field $V$ of the form (\ref{eq:V}) is to simply be continuous with compact suport in $\fh$; as for the derivatives $ZZp$ it suffices to consider them in their distributional sense.

\subsubsection*{Extension of quasiconformal deformations. \index{quasiconformal deformations!extension}}
The measurable Riemann mapping theorem in the Heisenberg group case is the infinitesimal analogue of the measurable Riemann mapping theorem of Ahlfors and Bers in the Heisenberg setting, but there is no result assuring the existence of a solution to the Beltrami system \index{Beltrami equation} of Equations (\ref{eq:B1}) and (\ref{eq:B2}). However, this theorem constitutes the key step to pass from quasiconformal deformations  \index{quasiconformal deformation} of the Heisenberg group $\fh$ to quasiconformal deformations  \index{quasiconformal deformation} of the complex hyperbolic plane $\bH^2_\C$. In the following we describe this passage, restricting ourselves to the smooth case. Assuming enough smoothness, quasiconformal mappings  \index{quasiconformal mapping}\index{mapping!quasiconformal} of $\bH^2_\C$ are necessarily symplectic transformations, i.e., diffeomorphisms $F$ such that $F^*\Omega=\Omega$, where $\Omega$ is the symplectic form in $\bH^2_\C$ derived by its K\"ahler metric. If $J$ is the 
natural complex structure in $\bH^2_\C$, then $F$ defines another complex structure $J_\mu=F^{-1}_*\circ J\circ F$ in $\bH^2_\C$ and there is an
associated complex antilinear self-mapping of the $(1,0)$-tangent bundle $T^{(1,0)}(\bH^2_\C)$ such that the $(1,0)$-tangent bundle of $J_\mu$ is $\{Z-\overline{\mu Z}\;|\;Z\in T^{(1,0)}\}$. The map $\mu$ is called the {\it complex dilation} of $F$ and there is a description of $\mu$ via a Beltrami system \index{Beltrami equation} of equations, see pp. 401--402 in \cite{DT}. The 
connection between quasiconformal symplectic transformations of the complex hyperbolic plane and quasiconformal contact transformations of the Heisenberg group is described as follows, see \cite{KR3} and \cite{KR4}:
 \begin{enumerate}
  \item [{(i)}] A (quasiconformal) symplectic transformation $F$ of the complex hyperbolic plane $\bH^2_\C$ extends to a (quasiconformal) contact transformation of the boundary.
  \item [{(ii)}] A quasiconformal deformation of the boundary extends to a quasiconformal deformation in the interior.
 \end{enumerate}
 In both cases, the constant of quasiconformality remains the same.

\section{Further Developments and Some Open Problems}\label{sec:further}

\subsubsection*{Extremal problems.}  In the classical theory, {\it extremal} quasiconformal mappings \index{extremal quasiconformal mapping}\index{quasiconformal mapping!extremal} are the ones minimising the maximal distortion \index{maximal distortion}(constant of quasiconformality) within a certain class of mappings in the complex plane or between Riemann surfaces. Since the times of Gr\"otzsch and Teichm\"uller, a method based on the modulus \index{modulus of families of curves} of curve families has been applied to detect  such mappings; it turned out that the very same method applies for  the mappings  which minimise a mean distortion functional in the class of quasiconformal mappings between annuli in the complex plane, \cite{BFP2}. Given a domain $\Omega$ in the complex plane,  a family $\mathcal{F}$ of mappings defined in $\Omega$ and a density $\rho$ corresponding to the geometry of $\Omega$, the {\it mean distortion functional} \index{distortion functional}\index{functional!distortion} $\mathfrak{M}_1(f,\rho)$ is 
$$
\mathfrak{M}_1(f,\rho)=\frac{\int_{\Omega} K(p,f)\rho(p)^2 d\calL^2}{\int_{\Omega}\rho(p)^2 d\calL^2}, \quad f\in\mathcal{F}.
$$
Recently in \cite{BFP1}, a variation of the modulus method  \index{modulus method} has been developed in the Heisenberg group setting by Balogh, \index{Balogh, Z.} F\"assler \index{F\"assler, K.} and Platis, \index{Platis, I.} to prove that there exists a minimiser of a mean distortion functional \index{distortion functional}\index{functional!distortion}  in the class of quasiconformal mappings \index{quasiconformal mapping}\index{mapping!quasiconformal} between Heisenberg spherical annuli. Given a domain $\Omega$ in the Heisenberg group, a family $\mathcal{F}$ of mappings defined in $\Omega$ and a density $\rho$ corresponding to the geometry of $\Omega$, the {\it mean distortion functional} $\mathfrak{M}_2(f,\rho)$ \index{distortion functional}\index{functional!distortion} is 
$$
\mathfrak{M}_2(f,\rho)=\frac{\int_{\Omega} K(p,f)^2\rho(p)^4 d\calL^3}{\int_{\Omega} \rho(p)^4 d\calL^3},\quad f\in\mathcal{F}.
$$
The minimiser in \cite{BFP1} is the {\it Heisenberg stretch map}, an extension of the usual stretch map \index{stretch map} 
\begin{equation}\label{eq:stretch}
f_k(z)=z|z|^{(1/K)-1},\quad K>1,
\end{equation}
of the plane. Moreover, the Heisenberg stretch map \index{Heisenberg stretch map} is essentially the unique minimiser of the mean distortion functional, \index{distortion functional}\index{functional!distortion} see \cite{BFP3}, a fact that is in strong contrast to the classical case where there exist infinite such minimisers, see \cite{Bel}. However, it does not minimise the maximal distortion  \index{maximal distortion}and this is also in contrast to the classical situation. The problem of finding such a minimiser is open. We note here that the modulus method  \index{modulus method} is up to now the unique tool for the detection of extremal mappings; in the Heisenberg setting results similar to Teichm\"uller's Existence and Uniqueness theorems are not available. Moreover, in the application of the modulus method \index{modulus method} there is an additional difficulty due to the lack of a Riemann Mapping Theorem. \index{Riemann Mapping Theorem} In the classical case, one can always reduce an extremal 
problem concerning mappings defined on a ring domain to an extremal problem about mappings defined in a circular annulus. But in the Heisenberg group case things are rather different: In \cite{KR6}, Kor\'anyi \index{Kor\'anyi, A.} and Reimann \index{Reimann, H. M.} calculated the capacity  \index{capacity}(that is, the modulus \index{modulus of families of curves} of the family of curves joining the two components of the boundary) of a ring whose boundary comprises two homocentric Kor\'anyi-Cygan spheres. \index{Kor\'anyi-Cygan sphere}\index{sphere!Kor\'anyi-Cygan} Rather 
surprisingly, a variety of rings centered at the origin (like rings whose boundary comprises of two Carnot-Carath\'eodory spheres), have exactly the same capacity, see \cite{Pla}. In general though, the calculation of moduli of curve families in an arbitrary ring inside $\fh$ is a difficult task.

Extremal problems also arise naturally in the theory of quasiconformal mappings  \index{quasiconformal mapping}\index{mapping!quasiconformal} of compact psudoconvex $\rm{CR}$ manifolds. \index{$\rm{CR}$ manifold}\index{manifold!$\rm{CR}$} Such mappings have been defined by Kor\'anyi \index{Kor\'anyi, A.} and Reimann  \index{Reimann, H. M.} in \cite{KR5} and the extremality problem can be stated as follows. Given two $\rm{CR}$ structures on a 3-dimensional contact manifold, determine the quasiconformal homeomorphisms that have the least maximal distortion  \index{maximal distortion}with respect to the two $\rm{CR}$ structures. Problems of this type have been studied by various authors, see for instance the works of Miner \cite{Mi}, \index{Miner, R.} and Tang \index{Tang, P.} \cite{T}.

\subsubsection*{Sobolev and H\"older exponents. \index{Sobolev exponent} \index{H\"older exponent}}
Besides being a minimiser of the mean distortion functional \index{distortion functional}\index{functional!distortion} and of the maximal distortion, the planar stretch map ({\ref{eq:stretch}) is extremal for various other problems.
Astala \index{Astala, K.} showed in \cite{A} that a planar $K$-quasiconformal mapping lies in the Sobolev space \index{Sobolev space} $W_{loc}^{1,p}$ with $p<{2K}/{(K-1)}$. The example of the stretch map demonstrates that the given bound is sharp and this proves Gehring's conjecture, \index{Gehring's conjecture}\index{conjecture!Gehring} see \cite{Geh}, on the exponent of higher integrability in the two-dimensional case.
Ahlfors showed in \cite{Ahl} that a planar $K$-quasiconformal mapping is locally H\"{o}lder continuous with exponent $1/K$; again the stretch map can be used to show that this exponent cannot be improved.

The analogue of Gehring's higher integrability result for quasiconformal mappings  \index{quasiconformal mapping}\index{mapping!quasiconformal} on the Heisenberg group has been established by Kor\'anyi-Reimann \index{Kor\'anyi-Reimann theory} in \cite{KR3}: A quasiconformal mappping $f$ on $\fh$ lies in $HW_{loc}^{1,p}(\Omega,\fh)$ for an exponent $p>4$. However, up to present there is no sharp upper bound known in the spirit of Astala. Using the Heisenberg stretch mapping, it is shown in \cite{BFP1} that 
$$
p(\fh,K)\leq \frac{4K^{\frac{1}{4}}}{K^{\frac{1}{4}}-1},
$$
where
\begin{equation*}
 p(\fh,K)=\sup\left\{p\geq 1:\; f\in HW_{loc}^{1,p}(\fh,\fh),\;f:\fh\to\fh\;K-\text{quasiconformal}\right\}.
\end{equation*}
On the H\"older continuity side, it is known (see p. 53 in \cite{KR3}) that quasiconformal mappings  \index{quasiconformal mapping}\index{mapping!quasiconformal} on the Heisenberg group  and on more general Carnot groups \index{Carnot group}\index{group!Carnot} (see \cite{H1}) are locally H\"{o}lder continuous. A bound for the H\"{o}lder exponent in terms of the distortion has been given in Theorem 6.6  in \cite{BHT} and it is not likely that we can use the Heisenberg stretch to prove that this bound is optimal. Therefore, what is the appropriate quasiconformal mapping?

\subsubsection*{Quasiconformal maps in abstract spaces.} The study of quasiconformal mappings  \index{quasiconformal mapping}\index{mapping!quasiconformal} in the Heisenberg group $\fh$ motivated the study of quasiconformal mappings  \index{quasiconformal mapping}\index{mapping!quasiconformal} to larger and  more abstract spaces; some of which are $\rm{CR}$ manifolds, metric  spaces with controlled geometry, \index{space with controlled geometry} Carnot groups \index{Carnot group}\index{group!Carnot} and sub-Riemennian \index{sub-Riemannian space} (Carnot-Carath\'eodory) spaces. It also gave rise to questions concerning the comparison between the classical Ahlfors-Bers and the Kor\'anyi-Reimann \index{Kor\'anyi-Reimann theory} theory. It is equally fascinating to detect the points where similarities do exist, but also the points where they break down; all these are briefly addressed below. We start from  spaces with controlled geometry, Loewner spaces \index{Loewner space} and Carnot groups. \index{Carnot 
group}\index{group!Carnot}
We refer the reader to the pioneering work of Heinonen \index{Heinonen, J.} and Koskela \index{Koskela, P.} \cite{HK1} and \cite{HK2}, as well as to the notes of Reimann \cite{R1}. In general, metric spaces with controlled geometry are the metric spaces which display some kind of regularity with respect to comparison of distance and volume; the latter is the essence of quasiconformal mappings.  \index{quasiconformal mapping}\index{mapping!quasiconformal} Such a metric space $(X,d)$ of dimension $Q>1$ (also known as an {\it Ahlfors-David regular metric space}) \index{Ahlfors-David space} is endowed with a Borel measure $\mu$ compatible with the metric $d$ in the following way: there exists a constant $C\ge 1$ such that for all metric balls $B_R$ of radius $R<{\rm diam}(X)$ the following inequality holds:
\begin{equation*}
 \frac{1}{C} R^Q\le \mu(B_R)\le C R^Q.
\end{equation*}
Quasiconformal mappings are defined in such spaces via the metric definition, \index{metric definition of quasiconformality} and the same holds for notions like the modulus \index{modulus of families of curves} of curve families and quasisymmetric mappings. There are  well defined notions of $Q$-modulus \index{modulus of families of curves} ${\rm Mod}_Q(\Gamma)$ of a family of curves $\Gamma$ and of quasisymmetry, entirely analogous to those holding in the Heisenberg setting. 

A metric space $(X,d)$ is called a $Q$-{\it Loewner space}, if there is a strictly increasing self-mapping $\eta$ of $(0,\infty)$ such that ${\rm Mod}_Q(\Gamma)\ge\eta(k)$, where $\Gamma$ is the family of curves connecting two continua $C_0$ and $C_1$ with
$$
\min\{{\rm diam}(C_0),{\rm diam}(C_1)\}\ge k{\rm dist}(C_0,C_1).
$$
The Heisenberg group $\fh$ endowed with the Carnot-Carath\'eodory metric \index{Carnot-Carath\'eodory metric} \index{metric!Carnot-Carath\'eodory} $d_{cc}$ is a 3-regular Loewner space, and a bigger class of Loewner spaces are the so-called {\it Carnot groups}. A Carnot group \index{Carnot group}\index{group!Carnot} is a simply connected nilpotent Lie group $N$ with a derivation $\alpha$ on its Lie algebra $\fn$ such that ${\rm ker}(\alpha)$ generates $\fn$. Via the exponential map, $N$ and subsequently $\fn$ are identified to $\mathbb{R}^m$ for some $m\in\mN$ and the group action is given by the Campbell-Hausdorff formula, see \cite{cap}.  The Haar measure of $N$ is just the Lebesgue measure of $R^m$ and a Carnot-Carath\'eodory distance \index{Carnot-Carath\'eodory distance}\index{distance!Carnot-Carath\'eodory} is well defined. The starting point of the study of Carnot groups \index{Carnot group}\index{group!Carnot} is probably Pansu's thesis \cite{Pa}, see also \cite{MM}. A quite extensive study of the 
topic of quasiconformal analysis of Carnot groups \index{Carnot group}\index{group!Carnot} is found in the work of Vodop'yanov, \index{Vodop'yanov, S.} see for instance \cite{V}, \cite{V2} and \cite{V3}. 

The primary problem that one is facing in the study of the above spaces is to give proper analytic and geometric definitions \index{geometric definition of quasiconformality} of quasiconformality which are equivalent to the general (and applying in all cases) metric definition. \index{metric definition of quasiconformality} In this  direction, see the works of Williams \cite{W}, \index{Williams, M.} Tyson \cite{T1} \index{Tyson, J.} and \cite{T2}, and Heinonen \index{Heinonen, J.} et al. \cite{HKST} On the other hand, the conditions of the metric definition itself can be substantially relaxed and this gives rise to quite striking results, see \cite{BKR} and the bibliography given there.

\subsubsection*{The holy grail.}
In contrast to the Teichm\"uller space \index{Teichm\"uller space} case where extremal quasiconformal mappings  \index{quasiconformal mapping}\index{mapping!quasiconformal} are used to describe the whole space, it seems that a lot of effort has to be made to  obtain (or not!) an analogous result for spaces like the complex hyperbolic quasi-Fuchsian space which is defined now. 
Complex hyperbolic quasi-Fuchsian space $\calQ_\C(\Sigma)$ of a closed surface $\Sigma$ of genus $g>1$ \index{complex hyperbolic quasi-Fuchsian space} is perhaps the most natural extension of the Teichm\"uller space \index{Teichm\"uller space} of $\Sigma$: it consists of representations of the fundamental group $\pi_1(\Sigma)$ into the isometry group ${\rm SU}(2,1)$ of complex hyperbolic plane which are discrete, faithful, totally loxodromic and geometrically finite. We underline here that those conditions prevent that space (as well as the real quasi-Fuchsian space, that is the space of discrete, faithful, totally loxodromic and geometrically finite representations of $\pi_1(\Sigma)$ into ${\rm PSL}(2,\C)$) to fall in the Mostow rigidity \index{Mostow's rigidity theorem.} setting; the representations are convex cocompact and not cocompact as in the assumptions of Mostow's rigidity theorem. \index{Mostow's rigidity theorem} In a convex cocompact representation the quotient of the convex hull of the limit set
has finite volume (so it may have infinite funnel ends but
no cusps). Thus the limit set can never be the
entire boundary; there is always a region of discontinuity. In particular, for quasi-Fuchsian or complex hyperbolic quasi-Fuchsian  groups
the limit set is a topological circle and there is a domain of
discontinuity in the boundary.

 There is a quite large bibliography on the subject. For a summary of results concerning $\calQ_\C(\Sigma)$ we refer the reader to \cite{PP1}. Perhaps the most prominent problem in the subject is to examine the analytical structure of $\calQ_\C(\Sigma)$. In the case of Teichm\"uller space \index{Teichm\"uller space} this is carried out via the Ahlfors-Bers theory and the challenge here is to use the Kor\'anyi-Reimann \index{Kor\'anyi-Reimann theory} theory of quasiconformal mappings  \index{quasiconformal mapping}\index{mapping!quasiconformal} in the Heisenberg group to obtain similar results. In this direction, and regardless of the lack of an existence theorem for the solution of the 
Beltrami equation, one is invited to start from an irreducible representation $\rho\in\calQ_\C(\Sigma)$ and to construct quasiconformal deformations  \index{quasiconformal deformation} on the Heisenberg group with starting point $\rho$, to determine exactly the tangent space at $\rho$ from the vector fields generating these deformations. The problem is still open (it has been named {\it the holy grail} by the researchers in the area).

\end{document}